\documentclass{article}

\usepackage{amsmath,amssymb,latexsym,amstext,amsthm}
\usepackage{dsfont}

\newtheorem{theorem}{Theorem}[section]
\newtheorem{corollary}[theorem]{Corollary}
\newtheorem{lemma}[theorem]{Lemma}
\numberwithin{equation}{section}

\newcommand{\CC}{\mathbb{C}}

\newcommand{\PP}{\mathbb{P}}
\newcommand{\ZZ}{\mathbb{Z}}
\newcommand{\NN}{\mathbb{N}}
\newcommand{\LL}{\mathcal{L}^{(\alpha,\beta)}}
\newcommand{\Lll}{\mathcal{\widetilde{L}}^{(\alpha,\beta)}}
\newcommand{\dgr}[1]{\mbox{{ \textrm{deg}\/}}[#1]}
\newcommand{\dsty}{\displaystyle}

\newcommand{\unifn}{\;\; {\mathop{\mbox{\LARGE $\rightrightarrows$}}_{n\to\infty}}\;\;}
\newcommand{\weak}{\;\;{\stackrel{\rm *}{\longrightarrow}}\;\;}

\title{Orthogonality  with respect to a Jacobi  differential operator and applications}

\author{Jorge Borrego Morell \thanks{Research partially supported by Dirección General de
Investigación, Ministerio de Ciencia e Innovación of Spain, under grant MTM2009-12740-C03-01}\\
Univ. Carlos III de Madrid, Spain\\ {E-mail: jbmorell@gmail.com} \\
\and Héctor Pijeira Cabrera \thanks{Research partially supported by  Ministerio de Econom\'{\i}a y Competitividad of Spain, under grant  MTM2012-36372-C03-01.}\\
Univ. Carlos III de Madrid, Spain \\ E-mail: hpijeira@math.uc3m.es}

\begin{document}

\maketitle

\begin{abstract}
Let  $\mu$ be a finite positive Borel measure on $[-1,1]$, $m$  a fixed natural number  and $\LL[f]=(1-x^2)f^{\prime \prime} + (\beta- \alpha -(\alpha+ \beta + 2)x) f^{\prime}$, with $\alpha, \beta > -1$. We study algebraic and  analytic  properties of the sequence of monic polynomials $(Q_n)_{n>m}$ that satisfy the orthogonality relations
$$\int \, \LL[Q_n] (x)x^k d\mu(x)=0 \quad
    \mbox{for all} \quad 0 \leq k \leq n-1.$$
    A fluid dynamics  model  for  source points location of a flow of an incompressible fluid with preassigned  stagnation points  is also considered.
\end{abstract}




\section{Introduction}

Let $\PP$ be the space of all polynomials. We say that a sequence $\{p_n\}_{n\in \ZZ_+} \subset \PP,$  is  a {polynomial system} if for each $n$ the $n$th polynomial $p_n$ is of degree $n$.

Let  $\mu$ be a finite positive Borel measure on $[-1,1]$  and  $\{L_n\}_{n\in \ZZ_+}$ the corresponding  system of monic orthogonal polynomials; i.e.
$L_n(z)= x^n+ \dots$ and
\begin{equation}\label{(1)}
 \langle L_n,L_k \rangle_{\mu}=\int L_n(x) L_k(x) d\mu(x) \left\{ \begin{array}{ll}
    \not =0 & \hbox{ if } n=k, \\
     = 0 & \hbox{ if }  n \neq k.  \\
   \end{array}\right.
\end{equation}

Denote by $\LL$ the Jacobi differential operator on the space $\PP$, with $\alpha, \beta > -1$, where
\begin{equation}\label{(2)}
\LL[f]=(1-x^2)f^{\prime \prime} + (\beta- \alpha -(\alpha+ \beta + 2)x) f^{\prime}, \quad f \in \PP,
    \end{equation}
or equivalently (cf. \cite[(4.2.2)]{Szg75})
 \begin{equation}\label{(3)}
\LL[f]=\frac{\left((1-x)^{\alpha+1}\,(1+x)^{\beta+1}\, f^{\prime}\right)^{\prime}}{(1-x)^{\alpha}\,(1+x)^{\beta}}, \quad f \in \PP.
    \end{equation}

From \eqref{(2)}  it follows  that $f$  and $\LL[f]$ are polynomials of the same degree. It is straightforward that integrating \eqref{(3)} with respect to the $(\alpha,\beta)$-Jacobi measure $d\mu_{\alpha,\beta}(x)= (1-x)^{\alpha}\,(1+x)^{\beta}\,dx$ \, on $[-1,1]$, we obtain
 \begin{equation}\label{(4)}
\int  \, \LL[f](x)\, d\mu_{\alpha,\beta}(x)=0, \quad f \in \PP.
    \end{equation}

We say that $Q_n$ is the $n$th monic orthogonal polynomial with respect to the pair  $(\LL, \mu)$  if
$\dgr {Q_n} \leq n$ and
\begin{equation}\label{(5)}
    \int \, \LL[Q_n] (x)x^k d\mu(x)=0 \quad
    \mbox{for all} \quad 0 \leq k \leq n-1.
\end{equation}

This nonstandard  orthogonality with respect to a differential operator was introduced in \cite{ApLoMa02}, where  the authors prove conditions of existence and uniqueness. The starting  points  of this work are \cite{PiBeUr10,BePiMaUr11}, where the orthogonality with respect to the differential operator $\mathcal{L}_{\zeta}[f]= f+ (z-\zeta) f^{\prime}$, $\zeta \in \CC$, was analyzed.

From \eqref{(1)}, we have that a monic polynomial $Q_n$ of degree $n$ is orthogonal with respect to  $(\LL, \mu)$
if and only if it is a polynomial solution of the differential equation
\begin{equation}\label{(6)} \LL[Q_n]=\lambda_n \,L_n, \; \mbox{ where }\;  \lambda_n=\lambda_{n}^{(\alpha,\beta)}=
    -n(1+n+\alpha+\beta).
\end{equation}

As will be shown, it is not always possible to guarantee the existence of a system of polynomials $\{Q_n\}_{n\in \ZZ_+}$ orthogonal with respect to the pair $(\LL,\mu)$.  Let   $m \in \NN$ be fixed, a fundamental role in this paper is played by the class   $\mathcal{P}_m(\alpha,\beta)$  defined as the finite positive Borel measures  $\mu$   supported on  $[-1,1]$  for which  there exist  a non negative  polynomial $\rho$ of degree $m$, such that $d\mu(x)= {\rho^{-1}(x)} d\mu_{\alpha, \beta}(x)$.

This manuscript deals with some algebraic and analytic aspects of the sequence of orthogonal polynomials $\{Q_n\}_{n>m}$ orthogonal with respect to the pair $(\LL,\mu)$. We provide asymptotic results for the sequence  of the orthogonal polynomials with respect to $(\LL,\mu)$ and study   the set of accumulation points of their zeros as well. In particular, we prove

\begin{theorem}\label{Th6}
Let  $\mu \in \mathcal{P}_m(\alpha,\beta)$, where $m \in \NN$. If $\{\zeta_n\}_{n>m}$ is a sequence of complex numbers  with limit  $\zeta \in \CC \setminus [-1,1]$ and $\{Q_n\}_{n>m}$ the sequence of monic orthogonal polynomials with respect to the pair  $(\LL, \mu)$ such that  $Q_n(\zeta_n)=0$, then the accumulation points of
zeros of $\{Q_{n}\}_{n>m}$ are located on the set $E=\mathcal{E}(\zeta) \bigcup [-1,1]$, where $\mathcal{E}(\zeta)$ is the ellipse
\begin{equation}\label{ZerosAcum}
    \mathcal{E}(\zeta):= \left\{ z \in \CC : z= \cosh(\eta_{\zeta}+i\theta), 0 \leq \theta  < 2 \pi\right\},
   \end{equation} and $\eta_{\zeta}:=\ln |\varphi(\zeta)|=\ln |\zeta+\sqrt{\zeta^2-1}|$.
    If $\dsty \delta(\zeta)=\inf_{-1\leq x \leq 1}|\zeta-x|>2$ then $E=\mathcal{E}(\zeta)$.
\end{theorem}

Let $\varphi(z)=z+\sqrt{z^2-1} $  be the function which maps the complement of $[-1,1]$ onto the exterior of the unit circle,  where we take the branch of $\sqrt{z^2-1}$ for which $|\varphi(z)|>1$ whenever $z \in \CC \setminus [-1,1]$. If $\mu \in \mathcal{P}_m(\alpha,\beta)$,  let $\nu_1, \nu_2, \ldots, \nu_m \in \CC$  be the $m$ zeros of the polynomial  $ \rho(z)= r \prod_{i=1}^{m} (z-\nu_i)$, for wich $d\mu(x)= \rho^{-1}(x) d\mu_{\alpha, \beta}(x)$.

For all $z \in \CC \setminus [-1,1]$ we define the function $\Phi(\rho, z) $ and a constant $\phi_m$ as
$$
    \Phi(\rho, z) =  \prod_{k=1}^{m} \frac{z-\nu_i}{\varphi(z)-\varphi(\nu_i)},  \quad  \phi_m =  2^m \; \exp\left( \frac{1}{2\pi} \int  \frac{\log(\rho(t))}{\sqrt{1-t^2}} dt\right) .
$$

The function $\phi_m \, \Phi(\rho, z)$ is a particular case of the well known \emph{Szeg\H{o}'s function} (cf.  \cite[\S 6.1]{Nev79}).

For the sequence of  monic orthogonal polynomials with respect to the pair $(\LL,\mu)$  the following asymptotic behavior holds

\begin{theorem}\label{RelativeAsymptotic}
Let $\{\zeta_n\}_{n>m}$  be a sequence of complex numbers  with limit  $\zeta \in \CC \setminus [-1,1]$, $m \in \NN$, $\mu \in \mathcal{P}_m(\alpha,\beta)$, and  $(Q_n)_{n>m}$ the sequence of monic orthogonal polynomials with respect to the pair $(\LL,\mu)$ such that $Q_n(\zeta_n)=0$, then:
\begin{enumerate}
  \item Uniformly on compact subsets  of  $\Omega=\{z \in \CC : |\varphi(z)| > |\varphi(\zeta)| \}$, (i.e. the exterior of the ellipse $\mathcal{E}(\zeta)$)
\begin{equation}
\label{AsintComp1}
\frac{Q_{n}(z)}{P^{(\alpha,\beta)}_n(z)}   \unifn   \phi^2_m \, \Phi(\rho, z).
\end{equation}
  \item Uniformly on compact subsets  of  $\Omega=\{z \in \CC : |\varphi(z)| < |\varphi(\zeta)| \} \setminus [-1,1]$
            \begin{equation}
\label{AsintComp2}
\frac{Q_{n}(z)}{P^{(\alpha,\beta)}_n(\zeta_n)}   \unifn   -\phi^2_m \, \Phi(\rho, \zeta).
\end{equation}
If $\dsty \delta(\zeta)>2$ then \eqref{AsintComp2} holds for $\Omega=\{z \in \CC : |\varphi(z)| < |\varphi(\zeta)| \}$ (i.e. the interior of the ellipse $\mathcal{E}(\zeta)$).
\end{enumerate}
\end{theorem}

The paper is organized as  follows. In Section \ref{AR} we study the existence of a system of polynomials $\{Q_n\}_{n>m}$, orthogonal with respect to the pair $(\LL,\mu)$.  Section \ref{sect3} is devoted to the study of recurrence relations and location of zeros of the polynomials $\{Q^{\prime}_n\}_{n>m}$. In   Sections \ref{sect4} and \ref{sect5}  we study the asymptotic behavior of the polynomials   $\{Q_n\}_{n>m}$, $\{\widehat{Q}_n\}_{n>m}$ and its zeros respectively. In the last section,  we introduce a fluid dynamics model for the interpretation of the critical points of $Q_n$.

\section{Existence and uniqueness}\label{AR}

It is  known that the  operator  $\LL$ has a {system of monic eigenpolynomials} $\{P^{(\alpha,\beta)}_n\}_{n\in \ZZ_+}$  and a sequence of constant $(\lambda_n)_{n\in \ZZ_+}$ (eigenvalues), such that
\begin{equation}\label{(7)}\LL[P^{(\alpha,\beta)}_n]=\lambda_n P^{(\alpha,\beta)}_n,
\end{equation}
where the eigenpolynomial $P^{(\alpha,\beta)}_n$ is  the $n$th monic orthogonal polynomial with respect to the $(\alpha,\beta)$-measure of Jacobi $d\mu_{\alpha,\beta}(x)= (1-x)^{\alpha}\,(1+x)^{\beta}\,dx$ on $[-1,1]$,  i.e.
\begin{equation}\label{(8)}
   \left\langle P^{(\alpha,\beta)}_n, x^k \right\rangle_{\alpha,\beta}= \int \, P^{(\alpha,\beta)}_n(x) \, x^k d\mu_{\alpha,\beta}(x)= 0, \mbox{ for all } 0 \leq k \leq n-1.
\end{equation}
Furthermore, from \cite[(4.21.6) and (4.3.3)]{Szg75}
\begin{eqnarray}
\nonumber \tau_n &=& \left\langle P^{(\alpha,\beta)}_{n}, P^{(\alpha,\beta)}_{n} \right\rangle_{\alpha,\beta} \\
\label{(24)}   &=& \frac{\Gamma(n+\alpha+1)\,\Gamma(n+\beta+1)\,\Gamma(n+1)\,\Gamma(n+\alpha+\beta+1)}{{2^{-(2n+\alpha+\beta+1)}}\,\Gamma(2n+\alpha+\beta+2) \, \Gamma(2n+\alpha+\beta+1)}.
\end{eqnarray}

\begin{theorem} \label{(Th1)} Let $n$ be a fixed natural number and  $\mu$   a finite positive Borel measure on $[-1,1]$. Then, the differential equation \eqref{(6)} has  a
monic polynomial solution $Q_n$ of degree $n$, which is unique up to an additive constant, if and only if
\begin{equation}\label{(11)}
\int L_n(x)d\mu_{\alpha,\beta}(x)=0,
\end{equation} where $L_n$ is as \eqref{(1)}.
\end{theorem}

\begin{proof} \

\noindent Suppose that there exists a polynomial $Q_n$ of degree $n$, such that $\LL[Q_n]=-n(1+n+\alpha+\beta) \,L_n$, where $L_n$ is the $n$th monic orthogonal polynomial for $\mu$. Then \eqref{(11)} is straightforward form \eqref{(4)}.

Conversely, suppose that  $L_n$, the $n$th monic orthogonal polynomial with respect to  $\mu$ satisfies  \eqref{(11)}. Let $Q_n$ be the polynomial of degree $n$ defined by
\begin{equation}
\label{(10)}
 Q_n(z) = P^{(\alpha,\beta)}_n(z) + \sum_{k=0}^{n-1}a_{n,k} P^{(\alpha,\beta)}_k(z),
\end{equation}
where $a_{n,0}$ is an arbitrary constant and
$$a_{n,k}= \frac{\lambda_{n}}{\lambda_{k}}\,b_{n,k}=\frac{\lambda_{n}}{\lambda_{k}\,\tau_k}\, {\left\langle L_n,
P^{(\alpha,\beta)}_k \right\rangle_{\alpha,\beta}}, \quad k=1,\ldots ,n-1.$$
From the linearity of $\LL[\cdot]$ and \eqref{(7)} we get that $\LL[Q_n]=-n(1+n+\alpha+\beta) \,L_n$.
\end{proof}

We are interested in discussing  systems of polynomials such that for all $n$ sufficiently large  they are solutions of \eqref{(6)}. In this sense, the next corollary is fundamental.

\begin{corollary} \label{(Cor1)} Let $\mu$ be a finite positive Borel measure supported on $[-1,1]$  such that   $d\mu(x)=r(x) d\mu_{\alpha,\beta}(x)$, where $r \in {L}^2(\mu_{\alpha,\beta})$. Then, $m$ is the smallest natural number such that for each $n>m$ there exists a monic polynomial $Q_n$ of degree $n$, unique up to an additive constant and orthogonal with respect to $(\LL, \mu)$, if and only if    $r^{-1}$ is a polynomial of degree   $m$.
\end{corollary}

\begin{proof} Suppose that $m$ is the smallest natural number such that  for each $n>m$ there exists a monic polynomial $Q_n$ of degree $n$, unique up to an additive constant and orthogonal with respect to $(\LL, \mu)$. According to  Theorem \ref{(Th1)}
$$\int  \frac{1}{r(x)} L_n(x)d\mu(x)=\int L_n(x)d\mu_{\alpha,\beta}(x)
\left\{ \begin{array}{cc}
 =0 & \mbox{if } n>m\\
 \neq 0 & \mbox{if } n=m
  \end{array}
\right..$$
But this is equivalent to saying that $\dsty
\frac{1}{r(x)}= \sum_{k=0}^{m} c_{k} L_k(x) $ with $c_m \neq 0$. The converse is straightforward.
\end{proof}

From the previous corollary,  if $\mu \in \mathcal{P}_m(\alpha,\beta)$ then the differential equation \eqref{(6)} has an unique monic polynomial solution $Q_n$ of degree $n$ for all $n>m$, except for an additive constant.

Let $\{\zeta_n\}_{n>m}$ be a sequence of complex numbers, where $m \in \NN$ is fixed,  and assume that $\mu \in \mathcal{P}_m(\alpha,\beta)$. We complement the definition of the sequence $\{Q_n\}_{n>m}$ in \eqref{(5)} by considering that henceforth $Q_n$ for each $n > m$ is  the polynomial  solution of the initial value problem
\begin{equation}\label{IVP_n}\left\{\begin{array}{rcl}
                              \LL[y] & = & \lambda_n \,L_n, \quad n > m,\\
                               y(\zeta_n) &=& 0.
                             \end{array}\right.
\end{equation}

We say that $\{Q_n\}_{n>m}$, is the sequence of monic orthogonal polynomials with respect to the pair  $(\LL, \mu)$ such that  $Q_n(\zeta_n)=0$.

If  $\widehat{Q}_n$ is  the monic polynomial of degree $n$ defined by the formula
\begin{equation}\label{(12)}
   \widehat{Q}_n(z)=\lambda_n\, \sum_{k=0}^{m} \frac{ b_{n,n-k}}{\lambda_{n-k}} \;  P^{(\alpha,\beta)}_{n-k}(z), \quad  b_{n,n-k}=\frac{1}{\tau_{n-k}} \, {\left\langle L_n,
P^{(\alpha,\beta)}_{n-k} \right\rangle_{\alpha,\beta}},
\end{equation}
then, the initial value problem \eqref{IVP_n} has the unique polynomial solution
\begin{equation}\label{(13)}
y(z)=Q_n(z)= \widehat{Q}_n(z) - \widehat{Q}_n(\zeta_n).
\end{equation}

\section{The polynomial $Q^{\prime}_{n}$}
\label{sect3}

Let $m \in \NN$, $\{\zeta_n\}_{n\in \NN}$ a sequence of complex numbers, and $\mu \in \mathcal{P}_m(\alpha,\beta)$ be fixed, then for all $n>m$ the polynomials $Q_n$ (solution of \eqref{IVP_n}) are uniquely determined by \eqref{(12)}--\eqref{(13)}. Without loss of generality, we will complete the sequence of polynomials $Q_n$ for all $n \in \NN$ as follows
\begin{eqnarray}
\nonumber {Q}_n(z) & = & \left(P^{(\alpha,\beta)}_{n}(z)- P^{(\alpha,\beta)}_{n}(\zeta_n)\right) \\ & & + \lambda_n\, \sum_{k=1}^{\min(m,n)} \frac{ b_{n,n-k}}{\lambda_{n-k}} \;  \left(P^{(\alpha,\beta)}_{n-k}(z)- P^{(\alpha,\beta)}_{n-k}(\zeta_n)\right), n\geq 1 \label{(17)}\\
\nonumber Q_0(z) & = &1
\end{eqnarray}
For convenience, only in  the previous formula we consider $\lambda_{0}=1$. Note that  $\{Q_n\}_{n\in \NN}$ defined by \eqref{(17)} is a system of polynomials, such that $Q_n(\zeta_n)=0$ for all $n\geq1$. Let us remark that if $n \leq m$, in general,  $\LL[Q_n] \neq \lambda_{n} L_n.$

Additionally, as  the degree of a polynomial is invariant  under the operator  $\LL[\cdot]$  and the polynomial $Q_n$,  for all $n \leq m$,  is of degree $n$ (see \eqref{(17)}),
\begin{equation}\label{(18)}
   (1, \LL \left[ Q_{1}\right], \ldots,  \LL \left[Q_{m}\right], L_{m+1}, \ldots, L_n, \ldots )
\end{equation}
is a  polynomial system.

In the following result,  we show that for  $n > (2m+1)$ the derivatives of the system of polynomials $Q_n$ satisfy  a  recurrence relation with a fixed  finite  number of terms.

\begin{theorem}\label{(Th2)}
Let $\mu \in \mathcal{P}_m(\alpha,\beta)$ where $m \in \NN$. Then
 if $R$ is any primitive of $\rho$, for each $n > (2m+1)$ the sequence of polynomials  $Q^{\prime}_{n}$  satisfies the  relation
\begin{equation}\label{(19)}
  R(z)Q^{\prime}_{n}(z)=\sum_{k=-m-1}^{m+1} \theta_{R,n,n-k}\, Q^{\prime}_{n-k}(z),
\end{equation}
where  the initial values $Q^{\prime}_{m+1}, \dots, Q^{\prime}_{2m+2}$ are given by the derivatives of \eqref{(12)} and
\begin{eqnarray}
\nonumber \theta_{R,n,n-k} & = & \frac{1}{\lambda_{n-k}}\,(\lambda_n\;  e_{R,n,n-k} + {d}_{n,n-k}),\\
 \label{(20)} e_{R,n,n-k} & = &  \frac{1}{l_{n-k}}\,{\dsty \int  \,R(x)\,L_{n-k}(x)\,L_{n}(x) \,d\mu(x)}, \\
\label{li-def}  l_i & = & \int   L^2_i(x)d\mu(x),
\end{eqnarray}
\begin{eqnarray}
\nonumber  d_{n,n-k} & = & \dsty \frac{1}{l_{n-k}}\sum_{j=j_1(k)}^{j_2(k)} \tau_{n-j}\, \widetilde{c}_{n,n-j} \, b_{n-k,n-j},\\
 \nonumber    & & j_1(k)=\max\{-1,k\} \mbox{ and } j_2(k)=\min\{m+1 ,m+k \}.\\
\nonumber    \widetilde{c}_{n,n-k} &=& \lambda_n\, \sum_{j=j_3(k)}^{j_4(k)}\, \frac{b_{n,n-j}\,{c}_{n-j,j-k}}{\lambda_{n-j}}, \\
\nonumber     & & j_3(k)=\max\{0,k-1\} \mbox{ and } j_4(k)=\min\{m,k+1\}, \\
\nonumber     c_{n,1}  &=&  -n,\\
\nonumber   c_{n,0}  &=&   \frac{2n\,(\alpha-\beta)\,(n+\alpha+\beta+1)}{(2n+\alpha+\beta)\,(2n+\alpha+\beta+2)},\\
\nonumber   c_{n,-1}  &=& 4n\, \frac{(n+\alpha)\,(n+\beta)\,(n+\alpha+\beta)\,(n+\alpha+\beta+1)}{(2n+\alpha+\beta)^2\,((2n+\alpha+\beta)^2-1)}.
\end{eqnarray}
\end{theorem}

Before starting the proof of Theorem \ref{(Th2)}, we will state and prove some lemmas.

\begin{lemma}
\label{(Lem1)}
Let $\mu \in \mathcal{P}_m(\alpha,\beta)$ where $m \in \NN$. Then for $n>m$
\begin{eqnarray}\label{(22)}
  L_n(z) &=&     \sum_{k=0}^{m}b_{n,n-k} \; P^{(\alpha,\beta)}_{n-k}(z), \\
 \label{(23)} \rho(z) P^{(\alpha,\beta)}_n(z) & = & \dsty  \tau_n \sum_{k=0}^{m} \frac{{b}_{n+k,n}}{l_{n+k}}  \; L_{n+k}(z),
\end{eqnarray}
where $b_{i,j} =  \frac{1}{\tau_j}\, \left\langle L_i, P^{(\alpha,\beta)}_j \right\rangle_{\alpha,\beta}$, $\tau_j$ as in \eqref{(24)} and $l_i$ as in \eqref{li-def}.
\end{lemma}

\begin{proof} As $\mu \in \mathcal{P}_m(\alpha,\beta)$ then  $\left\langle L_n, x^k  \right\rangle_{\alpha,\beta}=0$ for all $n>m$ and $k < n-m$.   Hence, $\dsty  b_{n,k} = 0$ for  $k=0,1,\ldots,n-m-1$ and \eqref{(22)} is established.

From the orthogonality relations \eqref{(1)} and \eqref{(8)}, if $i<j$ or $i>j+m$, we have
\begin{equation}\label{(26)}
   \int   L_i(x)  P^{(\alpha,\beta)}_{j}(x) d\mu_{\alpha,\beta}(x) =  \int   L_i(x)  P^{(\alpha,\beta)}_{j}(x) \, \rho(x) \, d\mu(x) =0.
\end{equation}
The relation \eqref{(23)} is  straightforward   from the Fourier expansion of $\rho \, P_n^{(\alpha,\beta)}$ in terms of the  $\left(L_k \right)$, $k=0,1,\ldots,n+m$ and \eqref{(26)}. \end{proof}

From the lemma above, the polynomials $\dsty \widehat{Q}_n$ defined in \eqref{(12)} and its derivatives can be written as a linear combination of the polynomials $\{L_n\}_{n\in \ZZ_+}$ as we show in the next lema.
\begin{lemma} \label{(Lem2)} Under the conditions of Lemma \ref{(Lem1)},  for $n>m$ the polynomials  $\dsty \widehat{Q}_n $  satisfy the following relations:
\begin{eqnarray}
 \label{(27)}     \int  \,  \widehat{Q}_n(x)  \, x^k \, d\mu_{\alpha,\beta}(x) & =& 0, \quad k=0,1,\ldots,n-m-1,\\
\label{(29)}(1-z^2) \rho(z) Q^{\prime}_n(z)  &= &  \dsty  \sum_{k=-m-1}^{m+1}\, d_{n-k,k}\;  L_{n-k}(z).
\end{eqnarray}
\end{lemma}

\begin{proof} From \eqref{(12)} the relations \eqref{(27)} and

\begin{equation}\label{(30)}
 (1-z^2)\,Q^{\prime}_n(z) =   (1-z^2)\,\widehat{Q}^{\prime}_n(z) =  \lambda_n\, \sum_{k=0}^{m} \frac{ b_{n,n-k}}{\lambda_{n-k}} \;  (1-z^2)\,\left( P^{(\alpha,\beta)}_{n-k}(z)\right)^{\prime},
\end{equation}
follows directly.

Using  the structure relation fulfilled by Jacobi polynomials (see \cite[(4.5.5)--(4.5.6)]{Szg75}), we have
$$
   (1-z^2)\left( P^{(\alpha,\beta)}_{n-k}(z)\right)^{\prime}=c_{n-k,1}\, P^{(\alpha,\beta)}_{n-k+1}(z) +c_{n-k,0}\, P^{(\alpha,\beta)}_{n-k}(z)+c_{n-k,-1}\, P^{(\alpha,\beta)}_{n-k-1}(z).
$$
Substituting this formula into \eqref{(30)}, we obtain
$$
    (1-z^2) Q^{\prime}_n(z) =  \sum_{k=-1}^{m+1}\, \widetilde{c}_{n,n-k}\;  P^{(\alpha,\beta)}_{n-k}(z),
$$
 and from \eqref{(23)}, \eqref{(29)} immediately follows. \end{proof}

\medskip

\begin{proof}[Proof of Theorem \ref{(Th2)}] As the sequence $\{Q_n\}_{n\in \ZZ_+}$  is a system of polynomials, then the sequence  of its derivatives $\{Q_n^{\prime}\}_{n\in \NN}$ is also system of polynomials.  Hence, the polynomial $R \, Q_{n}^{\prime}$ can be expanded as linear combination of the polynomials $\dsty \{\widehat{Q}^{\prime}_n\}_{n\in \NN}$, i.e. there exist  $(n+m)$  constants $\theta_{R,n,1}, \dots, \theta_{R,n,n+m}$  such that
\begin{equation}\label{(31)}
R(z)Q_{n}^{\prime}(z)=\sum_{k=-m}^{n-1}\,\theta_{R,n,n-k}\,\widehat{Q}^{\prime}_{n-k}(z).
\end{equation}

Let  $\Lll$ be the linear differential  operator on the space of all polynomials $\PP$  defined by  $\Lll[f^{\prime}]=\LL[f]$ for all $f \in \PP$, i.e.
$$
\Lll[f]=(1-x^2)f^{\prime} + (\beta- \alpha -(\alpha+ \beta + 2)x) f.$$

Since $\{L_n\}_{n\in \ZZ_+}$ is a system of polynomials and  $\Lll[\cdot]$  is a linear operator, the polynomial $\Lll[RQ_{n}^{\prime}]$ can be written as a linear combination of the system of polynomials \eqref{(18)} as follows
\begin{eqnarray}
\nonumber \Lll \left[R Q_{n}^{\prime} \right](z) &=& \sum_{k=-m-1}^{n-m-1}\,\theta_{R,n,n-k}\,\lambda_{n-k} \;L_{n-k}(z) \\ \label{(32)}
 & & +\sum_{k=n-m}^{n-1}\,\theta_{R,n,n-k} \;\LL \left[ Q_{n-k}(z) \right],
\end{eqnarray}

Taking $\Lll[\cdot]$ on the left hand side of the equality \eqref{(31)}, we get
\begin{eqnarray}
\nonumber \Lll[R Q_{n}^{\prime}](z)& = & R(z) \Lll[Q_{n}^{\prime}(z)]+(1-z^2) \rho(z)Q_n^{\prime}(z) \\
\label{(33)} & = & \lambda_n R(z) \; L_{n}(z)+(1-z^2) \rho(z)Q_n^{\prime}(z).
\end{eqnarray}
From \eqref{(1)}
\begin{eqnarray}\label{(34)}
    R(z) \; L_{n}(z) & = & \sum_{k=-m-1}^{m+1} \; e_{R,n,n-k}  \; L_{n-k}(z).
\end{eqnarray}
Substituting  \eqref{(29)} and \eqref{(34)} in \eqref{(33)}, we have
\begin{equation}\label{(35)}
\Lll[R Q_{n}^{\prime}](z) =  \sum_{k=-m-1}^{m+1} \; \left( \lambda_n\;  e_{R,n,n-k} + {d}_{n,n-k} \right)\; L_{n-k}(z).
\end{equation}

As $n \geq 2(m+1)$, we can assume that relation \eqref{(35)}  is the expansion of the polynomial $\Lll[R \, Q_{n}^{\prime}]$ in terms of the polynomials $L_n$. Now, identifying coefficients between  \eqref{(32)} and \eqref{(35)} we have that $\theta_{R,n,n-k}=0$ for all $k=1, \ldots, n-m-2 $, and we have the formulas \eqref{(19)}--\eqref{(20)}.
 \end{proof}

\section{Asymptotic behavior of the  sequence $\{\widehat{Q}_n\}_{n>m}$ and their zeros}
\label{sect4}

In this section we study the asymptotic behavior of the polynomials $\widehat{Q}_n$ and their zeros. Let us denote by $\| \cdot \|_{\alpha,\beta}$ and $\| \cdot \|_{\mu}$ the standard norms on  the spaces ${L}^2(\mu_{\alpha,\beta})$ and ${L}^2(\mu)$ respectively. The following  result  is   essential   in the proof of the theorems in  this and  the  next section.

\begin{theorem} \label{Th4}  Let $\mu \in \mathcal{P}_m(\alpha,\beta)$ where $m \in \NN$. Then
\begin{equation} \label{Q_Techo_Asym}
\frac{\widehat{Q}_n(z)}{P^{(\alpha,\beta)}_{n}(z)} \unifn \phi^2_m \, \Phi(\rho, z),
\end{equation}
 uniformly on closed subsets of $\overline{\CC} \setminus [-1,1]$  where $\overline{\CC}= \CC \bigcup \{\infty\}$.
\end{theorem}

\medskip First, we state a preliminary lemma  which  follows from   Theorems 26 and 29 of    \cite[\S 6.1]{Nev79}.
\begin{lemma} \label{LemmaAComp} Let $\mu \in \mathcal{P}_m(\alpha,\beta)$ where $m \in \NN$. If $\{L_n\}_{n\in \ZZ_+}$ is the sequence of monic orthogonal polynomials with respect to $\mu$
$$
    \frac{L_n(z )}{P^{(\alpha,\beta)}_{n}(z )} \unifn \phi^2_m \, \Phi(\rho, z),
$$
uniformly on closed subsets of  $\overline{\CC} \setminus [-1,1]$.
\end{lemma}

\begin{proof}[Proof of Theorem  \ref{Th4} ]
From  \eqref{(12)} and \eqref{(22)}
\begin{equation}\label{Relde Diferencia}
\frac{\widehat{Q}_n(z)-L_n(z)}{P^{(\alpha,\beta)}_{n}(z)}= \sum_{k=1}^{m} \left(\frac{ \lambda_n}{\lambda_{n-k}}-1\right)\,b_{n,n-k} \;  \frac{P^{(\alpha,\beta)}_{n-k}(z)}{P^{(\alpha,\beta)}_{n}(z)}.
\end{equation}

\noindent As  $\lambda_n= - n(1+n+\alpha+\beta),$  for each $k$ fixed, $k=1,2,\ldots,m$,
\begin{equation}\label{Opequenha}
\lim_{n \to \infty}\frac{ \lambda_n}{\lambda_{n-k}}= 1.
\end{equation}

Let $K$ be a closed subset of $\overline{\CC} \setminus [-1,1]$.  From the interlacing property of the zeros of consecutive Jacobi polynomials on $[-1,1]$, it easily follows that there exists a  constant $ M_*$  such that for all $z \in K$
\begin{equation}\label{Acotacion UniforJacobi}
 \left|\frac{P^{(\alpha,\beta)}_{n-k}(z)}{P^{(\alpha,\beta)}_{n}(z)} \right| <  M_k \leq M_*, \quad k=1,\ldots ,m,
\end{equation}
where $\dsty M_k=\sup_{\substack{ z \in K \\ x \in [-1,1]}} |z-x|^{-k}$ and $M_* = \max\{M_1,M_m\}$.

From  \eqref{(24)}, it is not hard to see that there exist two monic polynomials of degree $4(m-k)$ in the variable $n$, $q^{(\alpha,\beta)}_{1,4(m-k)}(n)$ and  $q^{(\alpha,\beta)}_{2,4(m-k)}(n)$, such that
$$  {\left\| P^{(\alpha,\beta)}_{n-k} \right\|^2_{\alpha,\beta}} = 4^{k-m} \frac{q^{(\alpha,\beta)}_{1,4(m-k)}(n)}{q^{(\alpha,\beta)}_{2,4(m-k)}(n)} \; {\left\| P^{(\alpha,\beta)}_{n-m} \right\|^2_{\alpha,\beta}}, \quad k= 1, 2,  \ldots, m.$$

Therefore, from the Cauchy-Bunyakovsky-Schwarz inequality and the extremal property of the monic orthogonal polynomials, for $n$ sufficiently large, we get
\begin{eqnarray}
\nonumber |b_{n,n-k}| &   \leq & \frac{\left\| L_n \right\|_{\alpha,\beta}}{\sqrt{\tau_{n-k}}}  \leq \sqrt{\frac{c_1}{\tau_{n-k}}} \, {\left\| L_n \right\|_{\mu}} \leq  \sqrt{\frac{c_1}{\tau_{n-k}}} \,  \frac{\left\| \rho \, P^{(\alpha,\beta)}_{n-m} \right\|_{\mu}}{|r|}
 \\  &  \leq &  \frac{c_1\, {\left\|  P^{(\alpha,\beta)}_{n-m} \right\|_{\alpha,\beta}} }{|r|\sqrt{\tau_{n-k}}}
 =  \frac{c_1 2^{m-k}}{|r|} \, \sqrt{\frac{q^{(\alpha,\beta)}_{2,4(m-k)}(n)}{ q^{(\alpha,\beta)}_{1,4(m-k)}(n)}}\leq \frac{{c_1}\,  2^{m+1}}{|r|}
\label{Acota Coef} \end{eqnarray}
where $ 1 \leq k \leq m$, $ \dsty c_1=\sup_{x \in [-1,1]} \rho(x)$ and $r$ is the leading coefficient of $\rho$. Hence by  \eqref{Relde Diferencia}, \eqref{Opequenha}, \eqref{Acotacion UniforJacobi}, and \eqref{Acota Coef}

\begin{equation*}
 \left|\frac{\widehat{Q}_n(z)-L_n(z)}{P^{(\alpha,\beta)}_{n}(z)} \right| \unifn 0, \mbox{ uniformly on closed subsets of }\, \overline{\CC} \setminus [-1,1],
\end{equation*}
and from Lemma \ref{LemmaAComp}  the asymptotic formula \eqref{Q_Techo_Asym} is established.
 \end{proof}

\begin{corollary}\label{Coro2} Let $\mu \in \mathcal{P}_m(\alpha,\beta)$ where $m \in \NN$. Then
\begin{enumerate}
  \item \begin{equation}\label{Q_techo_Raiz}
    \dsty \lim_{n  \rightarrow
\infty}\left|\widehat{Q}_{n}(z)\right|^{\frac{1}{n}}\,= \frac{|z+\sqrt{z^2-1}|}{2},
\end{equation}
uniformly on compact subsets of ${\overline{\CC}} \setminus [-1,1]$.
  \item The set of accumulation points of the zeros of the sequence of polynomials $\{\widehat{Q}_n\}_{n>m}$    is $[-1,1]$, i.e.
$$ \bigcap_{n\geq m} \overline{\bigcup_{k=n}^{\infty} \left\{ z: \widehat{Q}_k(z)=0 \right\}}=[-1,1].$$
For each $n$ at least $(n-m)$ zeros of $\widehat{Q}_n$ are contained on $[-1,1]$.
\end{enumerate}
\end{corollary}

\begin{proof} The first part of the theorem is an immediate consequence of relations \eqref{Q_Techo_Asym} and \cite[(8.21.9) and (4.21.6)]{Szg75}.

To prove the second part, from  \eqref{(27)} it easily follows that $\widehat{Q}_n$ has at least $n-m$ zeros   contained in $[-1,1]$  (cf. \cite[\S 3.3]{Szg75}).

The function in the right-hand  side of \eqref{Q_Techo_Asym} in Theorem \ref{Th4}  is holomorphic  and  does not have zeros in $\overline{\CC} \setminus [-1,1]$. Let $K$ be a closed subset of  $\overline{\CC} \setminus [-1,1]$, from the Rouch\'{e}'s theorem we have that  for $n$ large  the polynomial $\widehat{Q}_n$  does not have zeros on $K$, i.e. the  zeros of the sequence of polynomials $\{\widehat{Q}_n\}_{n>m}$    can not accumulate outside $[-1,1]$.

On the other hand, \eqref{Q_techo_Raiz}  implies the weak star asymptotic of the zero counting measures of the polynomials $\{\widehat{Q}_n\}_{n>m}$ (cf. \cite[Ch. 2]{StaTot92}). That is, if we associate to each $\widehat{Q}_{n}$ the measure
$
\dsty \mu_n=\frac{1}{n}\sum_{Q_{n}(\omega)=0}\delta_{\omega} \, ,
$ then
$ \dsty
d\mu_n(x)\weak \frac{1}{\pi}\frac{dx}{\sqrt{1-x^2}}
$
(the equilibrium distribution on $[-1,1]$) in the weak-* topology and this implies that the zeros of $\{\widehat{Q}_n\}_{n>m}$ must be dense in $[-1,1]$.  \end{proof}

\section{Asymptotic behavior of the  sequence $\{{Q}_n\}_{n>m}$ and their zeros}
\label{sect5}

Some basic properties of the zeros of  $Q_n$   follow  from \eqref{(6)}. For example, the multiplicity of the zeros of $Q_n$ is at most $3$, a zero of multiplicity $3$ is also a zero of $L_n$ and a zero of multiplicity $2$ is   a critical point of $\widehat{Q}_n$.

From the second part of Corollary \ref{Coro2},   we get that   $Q^{\prime}_n$ has at least $(n-m-1)$ zeros of odd multiplicity  on the open interval $]-1,1[$. For $m=1$ we have that

\begin{theorem} \label{Th3J} Under the same hypothesis of Theorem \ref{(Th2)}, if  $m=1$  the critical points of $Q_n$   interlace    the zeros of $L_n$.
\end{theorem}

\begin{proof}
  \quad  If $m=1$ then from \eqref{(27)} the polynomial $\widehat{Q}_n$ has at least $(n-1)$ real zeros of odd order on $]-1,1[$. But, $\widehat{Q}_n$ is a polynomial with real coefficients and degree $n$, consequently the zeros of $\widehat{Q}_n$ are real and simples. As ${Q}^{'}_n = \widehat{Q}^{'}_n$, from   Rolle's theorem all the critical points of ${Q}_n$ are real, simple  and $(n-2)$ of them are contained on $]-1,1[$.

  Denote  $\dsty P(x)= (1-x)^{\alpha+1}(1+x)^{\beta+1} Q_n^{\prime}(x)$. As $\alpha$ and $\beta$ are real numbers in general $P$ is not a polynomial.  Notice that $P$ is  a real--valued   differentiable  function on $[-1,1]$. Without loss of generality, suppose that there exists $x \in ]1,\infty[$ such that $P(x)=0$, as $P(1)=0$ from the Rolle's theorem there exists $x^{\prime} \in ]1,x[$ such that $P^{\prime}(x^{\prime})=0$. But, from \eqref{(3)} and \eqref{(6)} $\; P^{\prime}(x)=\lambda_n (1-x)^{\alpha}(1+x)^{\beta} L_n(x)$ and all the critical points of $P$ are contained in [-1,1]. Hence all the critical points of ${Q}_n$ are contained in $]-1,1[$.
  Again using the Rolle's theorem, it is straightforward that the  critical points of $Q_n$   interlace   the zeros of $L_n$.\end{proof}

From Corollary \ref{Coro2} we have that the set of accumulation points of $Q_n^{\prime}$ is $[-1,1]$. For $m=1$, Theorem \ref{Th3J} gives that the critical points of $Q_n$ are in $[-1,1]$, interlace   the zeros of $L_n$, and are simple. Numerical experiments also show this behavior for $m>1$. We conjecture that this always is the case.

For the proof of Theorems \ref{Th6} and  \ref{RelativeAsymptotic} we will use the following  result.

\begin{lemma} \label{Th5} Let  $\mu \in \mathcal{P}_m(\alpha,\beta)$, where $m \in \NN$ and define for $z \in \CC$,  $\dsty \Delta(z)= \sup_{x \in [-1,1]}|z-x|$ and
 $\dsty \delta(z)= \inf_{x \in [-1,1]}|z-x|$.. If $\{\zeta_n\}_{n>m}$ is a sequence of complex number  with limit  $\zeta \in \CC$ and $\{Q_n\}_{n>m}$ the sequence of monic orthogonal polynomials  with respect to the pair  $(\LL, \mu)$ such that  $Q_n(\zeta_n)=0$, then:
\begin{enumerate}
  \item For every  $\;d>1$ there is a positive number $N_d$, such that $\{z \in \CC : Q_n(z)=0 \} \subset \{z \in \CC : |z| \leq \Delta(\zeta)+d\}$ whenever $n >N_d$.
  \item If $\dsty \delta(\zeta) > 2$,  the zeros of $Q_n$ can not accumulate on $[-1,1]$ and for $n$ sufficiently large they are simple.
\end{enumerate}
\end{lemma}
\begin{proof}\
We already know that $Q_n(\zeta_n)=0$ and  if  $Q_n(z)=0$ then  $\widehat{Q}_n(z)=\widehat{Q}_n(\zeta_n)$. From the Gauss--Lucas theorem (cf. \cite[\S 2.1.3]{She02}), it is known that the critical points of $\widehat{Q}_n$ lie in the convex hull of its zeros and from  2. of Corollary  \ref{Coro2}  the zeros of the polynomials $\{\widehat{Q}_n\}_{n>m}$   accumulate on $[-1,1]$. Hence from the \emph{bisector theorem} (see the proof of Theorem  \ref{Th3J} or \cite[\S 5.5.7]{She02} )  $|z| \leq \Delta(\zeta_n)+1$ and the first  part of the theorem is established.

 To verify the second assertion of the theorem, note that if $z$ is a zero of $Q_n$, from  \eqref{(13)} we get
\begin{equation}\label{Cero_cond}
\prod_{k=1}^{n} \left|\frac{z-\widehat{x}_{n,k}}{\zeta_n-\widehat{x}_{n,k}} \right| = 1.
\end{equation}

Let $\dsty \mathcal{V}_{\varepsilon}([-1,1])= \{ z \in \CC : \delta(z) < \epsilon  \}$ be an $\varepsilon$--neighborhood  of $[-1,1]$.
On the other hand, as $\dsty \lim_{n \to \infty} \zeta_n=\zeta$, then for all $\varepsilon >0$ there is a $N_{\varepsilon}>0$ such that $|\delta(\zeta_n)-\delta(\zeta)|< \varepsilon $ whenever $n > N_{\varepsilon}$.

If  $\delta(\zeta)>2$, let us choose $\dsty \varepsilon = \varepsilon_{\zeta} = \frac{1}{2} \, \left(\delta(\zeta)-2 \right)$ and suppose that  there is a $z_0 \in \mathcal{V}_{\varepsilon_{\zeta}}([-1,1])$  such that $Q_n(z_0)=0$ for some $n > N_{\varepsilon_{\zeta}}$. Hence
\begin{equation}\label{Normal}
\prod_{k=1}^{n}\left|\frac{z_0-\widehat{x}_{n,k}}{\zeta_n-\widehat{x}_{n,k}} \right| < \left( \frac{2+\varepsilon_{\zeta}}{\delta(\zeta_n)} \right)^{n}< 1,
\end{equation}
 which is in contradiction with \eqref{Cero_cond}. Hence $\{z \in \CC : Q_n(z)=0 \} \bigcap \mathcal{V}_{\varepsilon_n}([-1,1])= \varnothing$ for all $n >N_{\varepsilon_{\zeta}}$, i.e.  the zeros of $Q_n$  can not accumulate on $\dsty \mathcal{V}_{\varepsilon_{\zeta}}([-1,1]).$

From \eqref{(13)} it is straightforward that a multiple zero of $Q_n$ is also a critical point of $\widehat{Q}_n$. But, from 2. of Corollary  \ref{Coro2} and the Gauss--Lucas theorem  the critical point of $\widehat{Q}_n$ accumulate on $[-1,1]$. Thus,  we have that for $n$ sufficiently large  the zeros of ${Q}_n$ are simple. \end{proof}

\begin{proof}[Proof of Theorem \ref{Th6}]
From    \eqref{(13)}  the  zeros of $Q_n$  satisfy
the equation
\begin{equation}\label{RaizCeros}
 \left|\widehat{Q}_{n}(z)\right|^{\frac{1}{n}}
 = \left|\widehat{Q}_{n}(\zeta_n)\right|^{\frac{1}{n}}\,.
\end{equation}
If $z \in \CC \setminus [-1,1]$, by taking limit when  $n \rightarrow \infty$, from 1. of lemma \ref{Th5}, and using  \eqref{Q_techo_Raiz} on  both sides of  \eqref{RaizCeros}, we have that the  zeros of the sequence of polynomials $({Q}_n)_{n>m}$    cannot accumulate outside  the set  $$\left\{ z \in \CC: |z+\sqrt{z^2-1}|=e^{\eta_{\zeta}}\right\} \, \bigcup \, [-1,1]. $$
Hence $z+\sqrt{z^2-1}= \, e^{\eta_{\zeta}+i\theta}$ and $z-\sqrt{z^2-1} = e^{-(\eta_{\zeta}+i\theta)}$ for $0 \leq \theta  < 2 \pi$, where we have that $2z = e^{\eta_{\zeta}+i\theta}+ e^{-(\eta_{\zeta}+i\theta)}.$

The assertion for  $\dsty \delta(\zeta)>2$ is straightforward from 2. of Lemma \ref{Th5}.
 \end{proof}

Now, we will  state the relative asymptotic between  the polynomials  $\{Q_n\}_{n>m}$  and the corresponding Jacobi polynomials $P^{(\alpha,\beta)}_n$.

\begin{proof}[Proof of Theorem \ref{RelativeAsymptotic}]\

\noindent 1.- Let us prove first that
\begin{equation}\label{AsintComp}
    \frac{Q_{n}(z)}{\widehat{Q}_n (z)} =1  - \frac{\widehat{Q}_n (\zeta_n)}{\widehat{Q}_n (z)} \unifn  1,
\end{equation}
uniformly on compact subsets $K$ of the set $\{z \in \CC : |\varphi(z)| > |\varphi(\zeta)| \}$.
In order to prove  \eqref{AsintComp} it is sufficient to show that
\begin{equation}\label{SufAsint}
\frac{\widehat{Q}_n (\zeta_n)}{\widehat{Q}_n (z)} \unifn 0, \quad \mbox{uniformly on } \; K.
\end{equation}
From  \cite[(8.21.9) and (4.21.6)]{Szg75}, we have the well known strong or power asymptotic of the monic Jacobi polynomials
\begin{equation}\label{StrongJacobiP}
 \frac{2^n\,P^{(\alpha,\beta)}_n(z)}{\varphi^n(z)} \unifn  \left(\frac{\varphi(z)-1}{2 (z-1)}\right)^{\alpha} \left(\frac{\varphi(z)+1}{2(z+1)}\right)^{\beta} \sqrt{\frac{\varphi^{\prime}(z)}{2}},
\end{equation}
 uniformly on  compact subsets  of  $\CC \setminus [-1,1]$. Note that
$$\frac{\widehat{Q}_n (\zeta_n)}{\widehat{Q}_n (z)}= \frac{\widehat{Q}_n(\zeta_n)}{P^{(\alpha,\beta)}_n(\zeta_n)}
\frac{P^{(\alpha,\beta)}_n(z)}{\widehat{Q}_n (z)}\frac{2^n\,P^{(\alpha,\beta)}_n(\zeta_n)}{\varphi^n(\zeta_n)} \frac{\varphi^n(z)}{2^n\,P^{(\alpha,\beta)}_n(z)} \left(\frac{\varphi(\zeta_n)}{\varphi(z)}\right)^n.$$
From \eqref{Q_Techo_Asym} and \eqref{StrongJacobiP} the first four  factors in the right hand side of the previous formula  have finite limits; meanwhile, the last factor tends to $0$  when $n \to \infty$, and we get \eqref{SufAsint}. Finally the assertion 1 is straightforward from \eqref{Q_Techo_Asym}.

\medskip \noindent 2.- For the assertion 2 of the theorem it is sufficient to prove that

 \begin{equation}\label{AsintComp-2}
    \frac{Q_{n}(z)}{\widehat{Q}_n (\zeta_n)} =\frac{\widehat{Q}_n (z)}{\widehat{Q}_n (\zeta_n)}  - 1 \unifn  -1,
\end{equation}
uniformly on compact subsets $K$ of the set $\{z \in \CC : |\varphi(z)| < |\varphi(\zeta)| \}\setminus [-1,1]$. Note that
$$\frac{\widehat{Q}_n (z)}{\widehat{Q}_n (\zeta_n)} = \frac{\widehat{Q}_n (z)}{P^{(\alpha,\beta)}_n(z)} \frac{P^{(\alpha,\beta)}_n(\zeta_n)}{\widehat{Q}_n(\zeta_n)} \frac{2^n\,P^{(\alpha,\beta)}_n(z)} {\varphi^n(z)}
\frac{\varphi^n(\zeta_n)}{2^n\,P^{(\alpha,\beta)}_n(\zeta_n)}
 \left(\frac{\varphi(z)}{\varphi(\zeta_n)}\right)^n.
$$
Now,  the first part of the assertion 2 is straightforward from \eqref{Q_Techo_Asym}.

If $\dsty \delta(\zeta)>2$,  let $\dsty \mathcal{V}_{\varepsilon}([-1,1])= \{ z \in \CC : \delta(z) < \epsilon  \}$ be a $\varepsilon$--neighborhood  of $[-1,1]$, where  $\dsty \varepsilon = \varepsilon_{\zeta} = \frac{\delta(\zeta)}{2}-1$. By the same reasoning that was deduced \eqref{Normal} we get that
\begin{equation}\label{Normal1}
\left|\frac{\widehat{Q}_n (z)}{\widehat{Q}_n (\zeta_n)}\right| < \kappa^{n}, \quad \mbox{for all } z \in \mathcal{V}_{\varepsilon}([-1,1]), \kappa <1.
\end{equation}
Hence from  the first part of the assertion 2 and  \eqref{Normal1} we get the second part of the assertion 2.
 \end{proof}

\section{Fluid dynamics model of sources and stagnation points}
\label{Sect2}

The fluid dynamic interpretation that we will consider in this section was introduced by H. Pijeira et al in \cite{BePiMaUr11}. In that  paper the hydrodynamic model was a reinterpretation  of the electrostatic model studied by H. Pijeira et al in \cite{PiBeUr10}.  The  difference between the fluid dynamic model  in  \cite{BePiMaUr11} and the model introduced in the present paper is  the complex potential   used.

Let us consider a flow of an incompressible fluid in the complex plane, due to a system of  $n-1$  source points ($n>1$) fixed at  $w_i$, $ 1 \leq i \leq n-1$, with unitary rate of fluid emission per unit time (\emph{strength of the source}),  and two additional source points at  $1$ and $-1$ with strength $a >0$ and $b >0$ respectively. Here, a \emph{source} is a point in which the fluid is continuously created and uniformly distributed in all directions with constant strength  (\emph{steady source}). Let us call \emph{flow field generated by a Jacobi set of sources} to a flow of a fluid  under the above conditions, or simple a \emph{flow field}.

The complex potential of a flow field  at any point $z$ (cf. \cite[Ch. 10]{Dur08} and \cite[Vol. II--Ch. 6]{Mark65}),  by the superposition principle of solutions, is given by
\begin{eqnarray}\nonumber
\Upsilon(z) &=& \sum_{i=1}^{n-1}
\log(z-w_i)+
a \, \log(z-1) +
b \, \log(z+1), \\  &=&  \log  \left( (z-1)^a \,(z+1)^b \, \prod_{i=1}^{n-1}(z-w_i) \right).\label{Complex_Pot}
\end{eqnarray}

From a complex potential $\Upsilon$, a \emph{complex velocity} $\mathcal{V}$ can be derived by differentiation ($\mathcal{V}(z)= \frac{d \Upsilon}{dz}(z)$). A standard problem associated with the complex velocity is to find the zeros, that correspond to the set of \emph{stagnation points}, i.e. points where the fluid has zero velocity.

We are interested in an inverse problem in the following sense, build a flow field such that the stagnation points are at preassigned points with \emph{nice} properties. As  it is well known, the zeros of orthogonal polynomials with respect to a finite positive Borel measure on $[-1,1]$ have a rich set of \emph{nice} properties (\cite[Chapter VI]{Szg75}),  and will  be taken as preassigned stagnation points. Here, we consider that $\mu \in \mathcal{P}_1(\alpha,\beta)$. In the next paragraph the statement of the problem will be established.

\textbf{Problem.} \emph{ Let $\{x_1,x_2, \ldots, x_n\}$ be the set of zeros of the $n$th orthogonal polynomial $L_n$ with respect to $\mu \in \mathcal{P}_1(\alpha,\beta)$ with $1<n$. Build a flow field (location of the  source points $w_1, \ldots, w_{n-1}$) such that the stagnation points  are attained at the points $ x_i$,   $i=1,2,\ldots,n$.}

Let $Q_{n}$ be a monic polynomial of degree $n$, and denote  its critical points by $\dsty \{w_1,w_2,\ldots,w_{n-1}\}$, thus
\begin{eqnarray*}
Q^{\prime}_{n}(z)  &=&  n\,\prod_{i=1}^{n-1}(z-w_i), \quad
\Upsilon(z) =   \log  \left(\frac{1}{n} (z-1)^a \,(z+1)^b\, Q^{\prime}_{n}(z)\right),\\
\mathcal{V}(z) &=& \frac{\partial \Upsilon}{\partial z}(z)  =   \frac{\left(  (z-1)^a \,(z+1)^b\, Q^{\prime}_{n}(z) \right)^{\prime}}{ (z-1)^a \,(z+1)^b\, Q^{\prime}_{n}(z)} =   \frac{\mathcal{L}^{(a-1,b-1)}[Q_{n}](z)}{ (z-1) \,(z+1)\, Q^{\prime}_{n}(z)}.
\end{eqnarray*}

From \eqref{Complex_Pot} and Theorem \ref{Th3J},  $ \dsty \frac{\partial \mathcal{V}}{\partial z}(x_k)=0$ for each stagnation point $x_k$ (zeros of $L_n$),  $k=1,2,\ldots,n$, i.e.
\begin{equation}\label{(15)}
 \mathcal{L}^{(a-1,b-1)}[Q_{n}](x_k)= 0, \quad
k=1,\,2,\, \cdots, n.
\end{equation}
From Corollary  \ref{(Cor1)}, there  exists a monic polynomial $Q_n$ of degree $n$, unique up to an additive constant,   satisfying    equation  \eqref{(15)},   i.e.
\begin{equation}\label{(16)}
 \mathcal{L}^{(a-1,b-1)}[Q_{n}](z)= \lambda_n L_n(z), \quad \lambda_n=-n(n+a+b-1).
\end{equation}

Note that \eqref{(16)} is  the same as \eqref{(6)} with $\alpha=a-1$ and $\beta=b-1$. Therefore the $n-1$ source point of the flow field   $\{w_1, \ldots, w_{n-1}\}$ are the critical point  of the $n$th orthogonal polynomial with respect to  the differential operator $\mathcal{L}^{(a-1,b-1)}$.

\textbf{Answer.} \emph{A flow fields generated by a Jacobi set of sources with complex potential \eqref{Complex_Pot} and preassigned stagnation points at  the zeros of the $n$th orthogonal polynomial with respect to the measure $\mu \in \mathcal{P}_1(\alpha,\beta)$ with  $n>1$, has its  sources points  (with unitary strength)    located at  the critical points of the $n$th orthogonal polynomial with respect to  $(\LL, \mu)$.}

In Theorem \ref{Th3J}, we proved that for  $m=1$ all the critical points of $Q_n$ are simple, contained in $[-1,1]$ and interlace  the zeros of $L_n$. At the  beginning of the Section \ref{sect5}, we conjectured that this theorem is true for all $m \in \NN$. If this were true, then it is not difficult to see that the above model holds for $m \in \NN$.

Note that, if we consider a system of electrostatic charges with potential given by \eqref{Complex_Pot}  instead of a system of   source points with the same potential function, then we have an analogous  electrostatic interpretation.

As it is known, the zeros of  the Jacobi polynomials have an electrostatic interpretation  (see \cite[\S 6.7]{Szg75}) as the equilibrium points of a certain potential function. For the case of orthogonality with respect to a differential operator  the electrostatic interpretation is an inverse problem in the sense that the equilibrium points are known and the question is to build the electrostatic field.

\section*{Acknowledgement.} The authors would like to thank Professor Guillermo L\'{o}pez Lagomasino  and the anonymous referees  for their careful revision of the manuscript and  suggestions which helped improve the presentation.

\end{document}